\documentclass[11pt]{amsart}

\usepackage{amsmath}
\usepackage{amssymb}
\usepackage{graphicx}
\usepackage{amscd}
\usepackage{stmaryrd}
\usepackage[mathscr]{euscript}

\theoremstyle{definition}

\numberwithin{equation}{section}

\newcommand{\sA}{\mathscr A}
\newcommand{\sB}{\mathscr B}
\newcommand{\sC}{\mathscr C}

\newcommand{\sE}{\mathscr E}

\newcommand{\sH}{\mathscr H}
\newcommand{\sJ}{\mathscr J}
\newcommand{\sK}{\mathscr K}

\newcommand{\sP}{\mathscr P}
\newcommand{\sR}{\mathscr R}

\title[Capacity and Quasicentral Modulus]{Capacity and the quasicentral modulus }

\author[D.-V. Voiculescu]{Dan-Virgil Voiculescu${}$}

\address{Department of Mathematics \\ University of California at Berkeley \\ Berkeley, CA\ \ 94720-3840}
\email{{\tt dvv@math.berkeley.edu}}

\begin{document}

\begin{abstract}
We point out that the quasicentral modulus is a noncommutative analogue of a nonlinear rearrangement invariant Sobolev condenser capacity. In the case of the shifts by the generators of a finitely generated group, the quasicentral modulus coincides with a corresponding nonlinear condenser capacity on the Cayley graph of the group.  Some other capacities  related to the quasicentral modulus are also discussed.
\end{abstract}

\maketitle


\section{Introduction}
\label{sec1}

The quasicentral modulus $k_{\sJ}(\tau)$ where $\tau$ is an $n$-tuple of operators and $(\sJ,|\ |_{\sJ})$ is a normed ideal is a number which is key in many questions about normed ideal perturbations. In particular $k_{\sJ}(\tau)$ is an essential ingredient in generalizations of the Kato-Rosenblum theorem and of the Weyl-von Neumann-Kuroda theorem to $n$-tuples of commuting hermitian operators. For a particular choice of the normed ideal ${\sJ}$ there are also connections with the Kolmogorov-Sinai entropy.

Recently studying the commutant ${\sE}(\tau;{\sJ})$ of $\tau$ mod ${\sJ}$, it turned out that this algebra provides new structure that helps understand the ubiquity of $k_{\sJ}(\tau)$.

Here we continue the exploration and we point out a noncommutative analogy with capacity in nonlinear potential theory. Note that $\sJ\cap{\sE}(\tau;{\sJ})$, which reduces to $\sJ$, but  can also be viewed as a noncommutative first order Sobolev space with respect to a rearrangement invariant norm $|\ |_{\sJ}$, while the role of the gradient goes to the $n$-tuple of inner derivations $[ \cdot , T_j]$ with respect to the components of $\tau$.

In the case of the regular representation of a finitely generated group and of the n-tuple of shifts in $\ell^2$ by the generators, the nonlinear potential theory analogy fits well with a connection to Yamasaki hyperbolicity and potential theory on the Cayley graph of the group which we had already noticed earlier.

\bigskip
\noindent
\underline{\qquad\qquad} 

2020 {\em Mathematics Subject Classification}. Primary: 46L89; Secondary: 31C45, 47L20.

{\em Key words and phrases}. quasicentral modulus, noncommutative nonlinear condenser capacity, normed ideals of operators.

\vfill
\noindent

\newpage

Besides the introduction and references there are six more sections. Section~2 recalls the definition of the quasicentral modulus. Section~3 provides background about condenser capacity in the nonlinear potential theory rearrangement invariant Sobolev space context. Then in section~4 we explain the noncommutative analogy. The analogy in the case of finitely generated groups is the subject of section~5.  Two capacities related to the quasicentral modulus are introduced in section~6. The brief concluding remarks in section~7 are about type $II_\infty$ factor and semifinite infinite von  Neumann algebra framework, a possible role for the ampliation homogeneity property and the Hilbert-Schmidt class and linear noncommutative potential theory.

This is a paper about heuristics. We take the first steps of developing the technical consequences of the analogy pointed out here in (\cite{25}).

\section{The quasicentral modulus}
\label{sec2}

Let $\sH$ be a  separable complex Hilbert space of infinite dimension and let $\sB(\sH)$, $\sK(\sH)$, $\sR(\sH)$ denote the bounded ,the compact and the finite rank operators on $\sH$
respectively.

If $(\sJ,|\ |_{\sJ})$ is a normed ideal (\cite{10}, \cite{15}) and $\tau = (T_j)_{1 \le j \le n}$ is a $n$-tuple of bounded operators on $\sH$, $k_{\sJ}(\tau)$ is defined as follows (\cite{18}, \cite{20}. \cite{21} or \cite{23}).                                                                                                                                                                                                                                                            

\bigskip

\begin{center}

\begin{minipage}{4in}
$k_{\sJ}(\tau)$ is the least $C \in [0,\infty]$, such that there exist finite rank operators $0 \le A_m \le I$ so that $A_m \uparrow I$ and we have
\[
\lim_{m \to \infty} \max_{1 \le j \le n} \vert [A_m,T_j] \vert_{\sJ} = C.
\]
\end{minipage}
\end{center}

Particularly important is the case of the ideals $\sC_p^-$ , $1 \le p \le \infty$ our notation of the Lorentz $(p,1)$ ideal where
\begin{center}
\[
\vert T \vert^-_p = \sum_{k\in{\mathbb N }} s_k k^{-1+1/p}
\]
\end{center}
with $s_1 \ge s_2 \ge \dots$ denoting  the eigenvalues of $(T^*T)^{1/2}$ in decreasing order. If $\sJ = \sC^-_p$ we denote the quasicentral modulus $k^-_p(\tau)$. Note that when $p = 1$ , $\sC_1 = \sC_1^-$. while if $p > 1$ the quasucentral modulus for the Schatten - von Neumann $p$-class can only take the values $0$ and $\infty$ (\cite {20}). In many cases sharp perturbation results are for $\sJ = \sC^-_p$, not for the Schatten - von Neumann class $\sC_p$ (\cite {20}, \cite{21}, \cite{23}).

\section {Condenser capacity background }
\label{sec3}
Let $\Omega \subset {\mathbb R}^n$  be an open set and let $K \Subset \Omega$ be a compact subset. The p-capacity of $K$ relative to $\Omega$ is  
\begin{center}
\[
cap_p (K; \Omega) = inf\{\int \mid \nabla u \mid ^p d \lambda \vert u \in C^\infty_0 (\Omega) , 0 \le u \le 1, u \vert_K \equiv 1 \}
\]
\end{center}
where $\nabla$ denotes the gradient, $d \lambda$ Lebesgue measure and $C^\infty_0 (\Omega)$ the $C^\infty$ functions with compact support in $\Omega$. The definition is then extended to open subsets $G \subset \Omega$ by taking the $sup$ over the compact subsets $K \Subset \Omega$ and then further to general sets $E \subset \Omega$ by taking the $inf$ over open sets $G \supset E$ (\cite {1}, \cite {9}, \cite {11}, \cite {13}). Note that a power - scaling, passing from $cap_p (K; \Omega)$ to $cap_p (K; \Omega)^{1/p}$ amounts to considering 
\begin{center}
\[
inf\{ \Vert \nabla u \Vert _p  \vert u \in C^\infty_0 (\Omega), 0 \leq u \leq 1, u \vert _K  \equiv 1 \}
\]
\end{center} 
where $\Vert \ \Vert_p$ is the $L^p$ norm on $\Omega$.

There is a natural generalization to the Lorentz-space setting, where $cap_{p,q}^{1/p} (K; \Omega)$ is defined as the
\begin{center}
\[
inf\{ \Vert \nabla u \Vert _{p,q}  \vert u \in C^\infty_0 (\Omega), 0 \leq u \leq 1, u \vert _K  \equiv 1 \}
\]
\end{center}
where $\Vert \ \Vert_{p,q}$ is the norm  of the Lorentz-space $L^{p,q} (\Omega)$  (\cite{6}, \cite{7}) and of course instead of the Lorentz $(p, q)$-norm, general rearrangement invariant norms can be considered. A natural variant of the definition of  $p$-capacity is with $u \in C^\infty_0$ replaced by $u \in W^{1,p}_0$ the closure of $C^\infty_0$ in the first order Sobolev space of functions with gradients in $L^p$. There are also similar variants of the definition of $cap_{p,q}(K:\Omega)$ involving Sobolev-Lorentz spaces (for more on Sobolev-Lorentz spaces (see \cite{4}, \cite{6}, \cite{7}, \cite{12}, \cite{17}).

We should also add that a condenser is a more general object involving two disjoint subsets of $\Omega$. Thus if $K, L \Subset \Omega$ are compact subsets and $K \cap L = \emptyset$ one considers 
\begin{center} 
\[
cap_p (K, L; \Omega) = inf\{\int \mid \nabla u \mid ^p d \lambda \vert u \in C^\infty_0 (\Omega), 0 \leq u \leq 1, u \vert _K  \equiv 1, u\vert_L \equiv 0 \}
\]
\end{center}
Similarly one defines $cap_{p,q} (K, L; \Omega)^{1/p}$  to be
\begin{center} 
\[                                                                                                                                                                                                                                                                                                           inf\{\Vert \nabla u \Vert _{p,q}  \vert u \in C^\infty_0 (\Omega), 0 \leq u \leq 1, u \vert _K  \equiv 1, u\vert_L \equiv 0 \}
\]
\end{center}
There are corresponding variants and extensions to the general rearrangement invariant setting like in the case $L = \emptyset$, which was our previous discussion.

\section{The noncommutative extension}
\label{sec4}

We shall first bring to the forefront the algebras of operators on $L^2(\Omega;d\lambda)$ in the definition of $cap_{p,q}(K, L: \Omega)^{1/p}$ or of a general Sobolev-rearrangement 
invariant analogue. Let $\sA = L^\infty (\Omega; d \lambda), \sB = C^\infty_0 (\Omega) , \sB \subset \sA $ which we identify with the algebras of multiplication operators on
 $L^2(\Omega; d\lambda)$ to which they give rise. Then $\sA$ is a von Neumann algebra and $\sB$ is a weakly dense *-subalgebra. For the sake of simplicity assume the two disjoint compact subsets $K, L \subset 
 \Omega$ are equal to the closures of their interiors. They then give rise to the idempotents $P = \chi_K, Q = \chi_L$ (the indicator functions of the sets). Thus we have $P = P^* = P^2, Q = Q^*= Q^2, PQ = 0$. Further let
\begin{center}
\[
\delta (f) = \Vert \nabla f \Vert_{p,q} \  {\text if} \  f \in C^\infty_0 (\Omega)
\]
\end{center}
or the same formula with some other rearrangement invariant norm instead of the $(p,q)$-Lorentz norm. Since $\nabla f$ has $n$ components we take $\max_{1 \le j \le n} \Vert \partial f/ \partial x_j \Vert_{p,q}$ or use some other norm on ${\mathbb C}^n$. Then  $\delta $ has the property 
\begin{center}
\[
\delta (fg) \le \Vert f \Vert \delta (g) + \Vert g \Vert \delta (f)
\]
\end{center}
that is, it is a differential seminorm (\cite{3}) or Leibniz seminorm (\cite{14}). In this simplified setting,the definition of  $(cap_{p,q} (K, L; \Omega))^{1/p}$ then becomes
\begin{center}
\[ 
inf\{\delta  (u) \vert u \in \sB, 0 \le u \le 1, P \le u  , uQ = 0 \}
\]
\end{center}.

The same construction can also be performed with noncommutative algebras. Let $\sH$ be a separable complex Hilbert space and let $\sA = \sB(\sH)$ and let $\sB = \sR(\sH)$. Let further $P, Q \in \sB$ be hermitian idempotents so that $PQ = 0$. If $\tau = (T_j)_{1 \le j \le n}$ is a $n$-tuple of operators in $\sB (\sH)$ we define the differential seminorm 
\begin{center}
\[
\delta (X) = \max_{1 \le j \le n} \vert [ T_j , X] \vert _{\sJ} 
\]
\end{center}
where $(\sJ , \vert \ \vert_\sJ )$ is a normed ideal. We get a quantity
\begin{center} 
\[
k_\sJ (\tau; P, Q) = inf \{ \delta (X) \vert 0 \le X \le I, X \in \sR (\sH), P \le X, XQ = 0 \}
\]
\end{center} 
which is a noncommutative analogue of the condenser capacity with respect to the Sobolev rearrangement invariant space. A detail we would like to point out is that since $ 0 \le X \le I $ in the definition, $P \le X$ is equivalent to XP = P. In case $ Q = 0 $ we get a quantity
\begin{center}
\[
k_\sJ (\tau; P) = k_\sJ (\tau; P, 0).
\]
\end{center}
Let $\sP (\sH) = \{ P = P^* = P^2 \vert P \in \sR (\sH) \}$ . It is clear that
\begin{center}
\[
P_1 \le P_2 \Rightarrow k_\sJ (\tau ; P_1) \le k_\sJ (\tau ; P_2)
\]
\end{center}
since passing from $P_1$ to $P_2$ decreases the set over which the $inf$ is taken. Finally the quasicentral modulus is recovered from these 
quantities as $k_\sJ (\tau ; I) $ that is
\begin{center} 
\[
k_\sJ (\tau) = sup \{ k_\sJ (\tau ; P) \vert P \in \sP (\sH) \} .
\]
\end{center}
Thus $k_\sJ  (\tau) $ is an analogue of a power-scaled Sobolev capacity of the whole ambient open set $\Omega$.

Obviously this construction can be performed with many other pairs $\sA, \sB$ and differential seminorms $\delta$, especially since we did not bother to put additional requirements on $\sA, \sB, \delta$ here. 

\section{Finitely generated groups}
\label{sec5}

Let $G$ be a group with a finite generator $\gamma = \{g_1,..., g_n \}$. Then $G$ is the set of vertices of its Cayley graph and $(g, g_jg)$ where $g \in G, 1 \le j \le n$ are the edges. Let further $(\sJ, \vert \ \vert_{\sJ} )$ be a normed ideal and let $\ell_{\sJ} (G)$ be the corresponding function space on $G$, that is functions $f : G \longrightarrow {\mathbb C}$ with the norm $\vert f \vert_\sJ$ being the $\sJ$-norm of the multiplication operator by $f$ in $\ell^2 (G)$. (That is the well-known bijection between normed ideals and symmetrically normed Banach sequence spaces.) By $r (G) \subset \ell^{\infty} (G)$ we denote the functions with finite support on $G$ and we consider the following differential seminorm on $r (G)$
 \begin{center} 
 \[
 \delta_\sJ (f) = \max_{1 \le j \le n} \vert f(g_j \cdot) - f(\cdot) \vert_\sJ .
\]
\end{center}
If $X \subset G$ is a finite subset, we define the $\sJ$-capacity of $X$ by
\begin{center}
\[
cap_\sJ (X) = inf \{ \delta_\sJ (f) \vert f \in r(G), 0 \le f \le 1, f \vert_X \equiv 1 \} .
\]
\end{center}
The definition is extended to general subsets $Y \subset G$ by
\begin{center}
\[
cap_\sJ (Y) = sup \{ \delta_\sJ (X) \vert X \subset Y, X \ {\text finite} \}.
\]
\end{center}
The $p$-hyperbolicity of Yamasaki (\cite{26}) of G, is the property 
\begin{center}
\[
cap_p (\{e\}) > 0
\]
\end{center}

where $e \in G$ is the neutral element and $cap_p$ is $cap_{\sC_p}$, where $\sC_p$  denotes the Schatten-von Neumann $p$-class. It is easily seen that $cap_\sJ (\{e\}) > 0$ is equivalent to 
\begin{center}
\[
X \ne \emptyset \Rightarrow cap_\sJ (X) > 0
\]
\end{center}
and this is equivalent to $cap_\sJ(G) > 0$. In particular the Yamasaki $p$-hyperbolicity is equivalent to $cap_p (G) > 0$. Similarly  Yamasaki $p$-parabolicity is the property $cap_p (\{e\}) = 0$ , which is equivalent to $cap_p (X) = 0$ for all $X \subset G$ and is also equivalent to $cap_p (G) = 0$.

The capacity construction has also a more general version for condensers

\begin{center}
\[
cap_\sJ (X_1 , X_2) = inf \{ \delta (f) \vert f \in r (G), 0 \le f \le 1, f \vert_{X_1} \equiv 1, f\vert_{X_2} \equiv 0 \}
\]
\end{center}

where $X_1 , X_2 \subset G$ are disjoint finite subsets.

In terms of algebras of operators these constructions correspond to $\sA = \ell^\infty (G), \sB = r (G) $ and the differential seminorm $\delta_\sJ$. Potential theory on graphs or networks being a well-developed subject ( see for instance \cite{2}, \cite{16}) it is natural to pass from Cayley graphs of groups to more general graphs, which also should apply to the above remarks. Of course in view of this kind of generalization the definition in the case of Cayley graphs should be modified by moving the $\max_{1 \le j \le n}$ in the definition of $\delta_\sJ $ so that the modified quantity does not depend on a labelling of the edges, for instance
\begin{center}
\[
\delta^*_\sJ (f) = \vert \max_{1 \le j \le n} \vert f(g_j \cdot) - f( \cdot ) \vert \vert_\sJ .
\]
\end{center}

On the other hand let $\lambda$ be the left regular representation of $G$ on $\ell^2 (G)$ and let $\lambda (\gamma)$ be the $n$-tuple of unitary operators $( \lambda (g_j))_{1 \le j \le n}$ . If $X \subset G$ is a finite subset, let $\chi_X \in r (G) $ be the indicator function, which being identified with the multiplication operator, is a hermitian projection in $\sR (\ell^2(G) ) \subset \sB (\ell^2(G)) $. We then have:
\begin{center}
\[
cap_\sJ (\{e\}) = k_\sJ (\lambda (\gamma) ; \chi_{\{e\}} )
\]
\end{center}

 \begin{center}
 \[
cap_\sJ (X) = k_\sJ ( \lambda (\gamma) ; \chi_X)
\]
\end{center}

\begin{center} 
\[
cap_\sJ (X_1, X_2)  = k_\sJ (\lambda (\gamma) ; \chi_{X_1} , \chi_{X_2} )
\]
\end{center}

\begin{center}
\[ 
cap_\sJ (G) = k_\sJ (\lambda (\gamma)) 
\]
\end{center}
where $X_1, X_2 \subset G$ are disjoint finite subsets. Clearly the LHS's are $\ge$ than the RHS's because they involve infimums over subsets of the sets over which we take infimums in the RHS's. For the converse one uses the projection of norm one 
\begin{center}
\[
diag : \sB(\ell^2 (G))  \rightarrow \ell^\infty (G)
\]
\end{center}
which amounts to taking the diagonal of the matrix of an operator with respect to the canonical basis and which decreases $\sJ$-norms. Remark further that
\begin{center}
\[
\vert [\lambda (g) , A] \vert_\sJ = \vert \lambda (g) A \lambda (g^{-1}) - A \vert_\sJ     
\]      
\end{center}
which is
\begin{center}
\[                                                                                                                                                                                                                              \ge \vert diag (\lambda (g) A \lambda (g^{-1}) - A ) \vert_\sJ = \vert \lambda (g) (diag A) \lambda (g^{-1}) - diag A  \vert_\sJ
\]
\end{center}

so that if the matrix of A has only finitely many non-zero entries one sees that these kind of quantities in the infimums of the RHS's majorize quantities in the infimums of the LHS's. Combining this with an approximation argument  involving  compressions  by projections corresponding to certain finite sets, of the finite rank operators one obtains operators with matrices with finitely many non-zero entries which completes the proof.  

\section{Capacities related to the quasicentral modulus}
\label{sec6}

In case $\tau$ is a $n$-tuple of commuting hermitian operators, there are two capacities, in the usual sense, that arise in connection with the quasicentral modulus. Here, it will be suitable to describe the norm $\vert \ \vert_\sJ$ by a norming function $\bf {\Phi}$, that is
\begin{center}
\[
\vert X \vert_\sJ = \bf {\Phi} (s_1, s_2, ... )
\]
\end{center}
where $ s_1 \ge s_2 \ge ... $ are the eigenvalues of $(X^* X)^{1/2} $ in decreasing order. 

Up to a power-scaling, the Hausdorff measure has a straightforward generalization involving a norming function. If $E \subset {\mathbb R}^n $ is a bounded set, let 
\begin{center}
\[
{\bf U}_{\bf {\Phi}} (E) = inf \{  {\bf  \Phi} (r_1, r_2, ...) \vert \bigcup_{j \in {\mathbb N}} B(x_j ; r_j) \supset E \}
\]
\end{center}
and if $\epsilon > 0$, let ${\bf U}_{{\bf {\Phi}},{\epsilon}}$ be the above infimum restricted to coverings with balls of radius $< \epsilon$. We then take the limit as $\epsilon \rightarrow 0$  of the ${\bf U}_{{\bf {\Phi}}, {\epsilon}}$ and we arrive at a quantity 
${\bf U}^*_{\bf {\Phi}} (E)$, which is a generalization of power scaled Hausdorff measure.

Let $\tau$ be a $n$-tuple of commuting hermitian operators with a cyclic vector $\xi$ and let $\sigma (\tau) \subset {\mathbb R}^n$ be its spectrum. Then ${\bf U}^*_{\bf {\Phi}} (\sigma (\tau))$ provides an upper bound for $k_{\bf {\Phi}} (\tau) $. We sketch below the argument, which is along the lines of \cite{18}, \cite{23}.

If $\sigma (\tau) \subset \bigcup_{j \in {\mathbb N}} B(x_j ; r_j) , r_j < \epsilon$, let $(\omega_j)_{j \in {\mathbb N}} $ be a partition of $\sigma (\tau)$ into Borel sets so that $\omega_j \subset B(x_j ; r_j)$. Let $P_n \in \sP (\sH)$ be the projection onto $\bigoplus_{1 \le j \le n} {\mathbb C} E (\omega_j) \xi $, where $E(\cdot)$ denotes the spectral measure of $\tau$. Then 
\begin{center}
\[
\vert [ \tau , P_n] \vert_{\bf {\Phi}} \le 2 {\bf {\Phi}} (r_1, ... , r_n, 0, ...)
\]
\end{center}
and
\begin{center}
\[
{\Vert [\tau, P_n] \Vert \le 2\epsilon},   {n \rightarrow \infty} \Rightarrow {\Vert P_n \xi - \xi \Vert \rightarrow 0} 
\]
\end{center}
(the norms for $n$-tuples are the $max$ of the norms of the components). Letting $ \epsilon \rightarrow 0$ we can find a sequence $Q_m \in \sP (\sH)$ so that
\begin{center}
\[
\Vert Q_m \xi - \xi \Vert \rightarrow 0 , \Vert [ Q_m , \tau ] \Vert \rightarrow 0 
\]
\end{center}
and 
\begin{center}
\[
\limsup_{m \rightarrow \infty} \vert [ Q_m , \tau] \vert_{\bf {\Phi}} \le 2 {\bf U}^*_{\bf {\Phi}} (\sigma (\tau)) .
\]
\end{center}
Note that if $Q$ is a weak limit of $Q_m$ 's we have $Q \xi = \xi, [ Q, \tau ] = 0$ which implies $Q = I $. This suffices to imply $k_{\bf {\Phi}} (\tau) \le 2 {\bf U}^*_{\bf {\Phi}} (\sigma (\tau))$ (\cite{18}).

In view of the results about $k_{\bf {\Phi}}$ in \cite{20}, the case of norming functions $\bf {\Phi}$ defined from a sequence $\pi_j \downarrow 0$ as $j \rightarrow \infty$, by the formula
\begin{center}
\[
 {\bf {\Phi}} (r_1, r_2, ... ) = \sum_j \pi_j r_j
 \]
 \end{center}
 where $r_1 \le r_2 \le ... $ is of particular interest.

 It is also possible to proceed in the opposite direction, from the quasicentral modulus to a set function on Borel sets.

 Let ${\mu}$ be a Radon measure on ${\mathbb R}^n$ with compact support and let $E \subset {\mathbb R}^n$ be a Borel set. On $L^2 (E, {\mu} \vert E) $ we consider $\tau_E$ the n-tuple of multiplication operators by the coordinate functions. We define 
 \begin{center}
 \[
{ {\bf Q}_{\mu , {\bf {\Phi}}}} (E) = k_{\bf {\Phi}} (\tau_E).
 \]
 \end{center}
 Note that ${\bf Q}_{\mu , {\bf {\Phi}}} (E)$ depends only on the absolute continuity class of $\mu$, that is, if $\mu$ and $\nu$ have the same 
 null sets, then
 \begin{center}
 \[
 {{\bf Q}_{\mu , {\bf {\Phi}}}} (E) = {{\bf Q}_{\nu , {\bf {\Phi}}}} (E).
 \]
 \end{center}
Particular instances of this construction were used in \cite{18}, \cite{23} to study $k^-_p$. If ${\bf {\Phi}}^-_p$ is the corresponding norming function (i.e. the case of the sequence $\pi_j = j^{-1 + 1/p}$ ), then if $p = n$ we proved in \cite{18} that 
\begin{center}
\[
{{\bf Q}_{\mu, {\bf {\Phi}}}} (E) = C( \lambda (E))^{1/n}
\]
\end{center}
when $\mu \vert E$ is equivalent to $n$-dimensional Lebesgue measure. We found a similar result in non-integer dimensions with Lebesgue measure replaced by the Hutchieson measure on certain selfsimilar fractals, on which it coincides with Hausdorff measure \cite{23}.

This idea of this construction can also be adapted to bounded metric measure spaces, but without some $n$-tuple of functions playing the role of the coordinate functions, it may be more natural to use the family of all real-valued Lipschitz functions with Lipschitz constant $\le 1$ and to consider $k_{\bf {\Phi}} $ of the infinite family of multiplication operators arising from these functions in $L^2 (E, {\mu} \vert E)$ where E is a Borel set.

Note also that to prove that a given ${\bf Q}_{\mu, {\bf {\Phi}}}$ is $> 0$, in view of \cite{20} requires more difficult analytic results. like in \cite{8}, \cite{18}.

\section{Concluding remarks}
\label{sec7}

\medskip
\noindent
{\bf A). Infinite semifinite von Neumann algebras, factors of type ${II}_{\infty}$}

It is obvious that the quasicentral modulus can be defined in type ${II}_\infty$ factors. One takes $\sA$ to be the factor, instead of $\sB (\sH)$ and $\sB$ the ideal of finite rank bounded operators with respect to the trace of the factor, that is operators in the factor which have support projections with finite trace. There are also normed ideals of such factors  and the differential seminorm is also defined as the $max$ of the ideal norms of the commutators with the components of the $n$-tuple of operators. The problem is to develop good examples which can be studied, part of the difficulties arising from the existence of non-isomorphic such factors, the multitude of non-conjugate masa's etc On the other hand, the definition of the quasicentral modulus can be further extended to separable semifinite infinite von Neumann algebras with a specified normal semifinite faithful trace and $n$-tuples of trace-preserving automorphisms $(\alpha_j)_{1 \le j \le n}$. The differential seminorm for a bounded finite rank element $A$ being
\begin{center}
\[
\delta (A) = \max_{1 \le j \le n} \vert \alpha_j (A) - A \vert_{\sJ}
\]
\end{center}
where $\vert \ \vert_{\sJ}$ is a corresponding normed ideal norm. Examples would include commutative von Neumann algebras $\sA$. This easily translates to $n$-tuples of measure-preserving automorphisms of measure spaces with an infinite $\sigma$-finite measure and a specified symmetric function space norm (assuming the measure has no atoms). Perhaps such invariants have the advantage of having many examples in sight when compared to the factor situation.

\medskip
\noindent
{\bf B). Power scaling}

The analogy of Sobolev capacities with respect to rearrangement invariant norms and the quasicentral modulus is up to a power scaling. One may wonder whether there are natural candidates of similar exponents for the quasicentral modulus. Introducing these  exponents may improve the analogy. For instance in the case of $k^-_p$, which corresponds to the Lorentz $(p,1)$ ideal $\sC^-_p$ one should consider $(k^-_p)^p $. The results of \cite{18}, \cite{19}, \cite{24} for commuting $n$-tuples of hermitian operators show certain advantages in the case of integer dimension and also in certain examples with non-integer dimension. These results suggest that in case the normed ideal $(\sJ , \vert \ \vert_{\sJ})$ has the ampliation homogeneity property \cite{24}  
\begin{center}
\[
k_{\sJ} (\tau \otimes I_m) = m^{1/{\alpha}} k_{\sJ} (\tau) 
\]
\end{center}
the exponent should  be $\alpha$.        

\medskip
\noindent
{\bf C). Commutants mod Hilbert-Schmidt}

In (\cite{22}, 6.3) we found that commutants mod $\sC_2$ of normal operators were objects which were related to the noncommutative potential theory based on Dirichlet forms (\cite{5}). One may see in this a certain similarity with what happens classically with $p$-capacity when $p$ = 2.



\end{document}